\documentclass[12pt]{article}
\usepackage[final]{epsfig}
\usepackage{graphics}
\usepackage{amsmath}
\usepackage{amsfonts}
\usepackage{latexsym}
\usepackage{amssymb}
\usepackage{graphicx}
\usepackage{url}
\usepackage{epstopdf}

\newtheorem{lemma}{Lemma}[section]
\newtheorem{proposition}[lemma]{Proposition}
\newtheorem{remark}[lemma]{Remark}

\newtheorem{theorem}{Theorem}

\newtheorem{corollary}[lemma]{Corollary}

\newcommand{\proofend}{$\Box$\bigskip}

\newcommand{\R}{{\mathbb R}}

\def\proof{\paragraph{Proof.}}
\def\Ker{\mathop{\rm Ker}}

\title{Iterating evolutes of spacial polygons and of spacial curves}

\author{Dmitry Fuchs\footnote{
Department of Mathematics, 
University of California, 
Davis, CA 95616;
 fuchs@math.ucdavis.edu}
 \and 
 Serge Tabachnikov\footnote{
Department of Mathematics,
Pennsylvania State University,
University Park, PA 16802;
tabachni@math.psu.edu}
} 

\date{}

\begin{document}

\maketitle

\begin{abstract}
The evolute of a smooth curve in an $m$-dimensional Euclidean space is the locus of centers of its osculating spheres, and the evolute of a spacial polygon is the polygon whose consecutive vertices are the centers of the spheres through the consecutive $(m+1)$-tuples of vertices of the original polygon. We study the iterations of these evolute transformations. This work continues the recent study of similar problems in dimension two, see \cite{ArFuIzTaTs}. Here is a sampler of our results. 

The set of $n$-gons with fixed directions of the sides, considered up to parallel translation, is an $(n-m)$-dimensional vector space, and the second evolute transformation is a linear map of this space. If $n=m+2$, then the second evolute is homothetic to the original polygon, and if $n=m+3$, then the first and the third evolutes are homothetic. In general, each non-zero eigenvalue of the second evolute map has even multiplicity. We also study curves, with cusps, in 3-dimensional Euclidean space and their evolutes. We provide continuous analogs of the results obtained for polygons, and present a class of curves which are homothetic to their second evolutes; these curves are spacial analogs of the classical hypocycloids.
\end{abstract}

\section{Introduction}\label{intro}

The evolute of a smooth plane curve is the locus of its centers of curvature or, equivalently, the envelope of the family of its normal lines. The construction of the evolute can be iterated. In our recent paper \cite{ArFuIzTaTs}, we studied such iterations; we also investigated discrete versions of this problem where smooth curves are replaced by polygons.

In the present paper we study higher dimensional analogs of this problem, the iterations of the evolutes of polygons in $\R^m$ and of  closed smooth curves, possibly, with cusps, in $\R^3$. Our investigation consists of two parts, concerning  polygons and curves, respectively. 

The first part of the paper concerns polygons in Euclidean spaces of arbitrary dimensions.

For an $n$-gon in $\R^m$ with $n \geq m+2$,  we define the evolute as the $n$-gon whose vertices are the centers of the spheres through the $(m+1)$-tuples of its consecutive vertices. (There exists an alternative definition of an evolute of an $n$-gon as an $n$-gon whose vertices are the centers of spheres tangent to the $(m+1)$-tuples of consecutive sides; we briefly consider this variant but do not study it in any detail.)

A polygon ${\bf P}$ is an involute of a polygon ${\bf Q}$ if ${\bf Q}$ is the evolute of ${\bf P}$. The existence and uniqueness of an involute of a generic $n$-gon in $\R^m$ depends on $n$ and $m$ (and, of course, on which of the two definitions of the evolute is chosen). In Proposition \ref{existinvol}, we give a complete answer to this question.

For an $n$-gon, its tangent indicatrix (the image of the tangent Gauss map) is a $n$-gon in $S^{m-1}$ or, equivalently, a cyclically ordered sequence of $n$ unit vectors in $\R^m$. Given a spherical $n$-gon ${\bf v}$, the space ${\cal P}_{\bf v}$ of the spacial $n$-gons with the respective directions of the sides, considered modulo parallel translations, is an $(n-m)$-dimensional vector space,  with the signed side lengths as coordinates.

The sides of the evolutes of the polygons in ${\cal P}_{\bf v}$  also have fixed directions, 
and the respective spherical polygon ${\bf u}$ is spherically dual to ${\bf v}$. The
sides of the second evolute  are parallel to the sides of the original polygon.  

The evolute map ${\cal P}_{\bf v} \to {\cal P}_{\bf u}$ is a linear map of $(n-m)$-dimensional vector spaces, depending on ${\bf v}$ (Theorem \ref{linear}). Generically, this map has full rank if $n-m$ is even, and has a 1-dimensional kernel if $n-m$ is odd.
The second evolute map  is a linear self-map of ${\cal P}_{\bf v}$. 

For a generic spherical $n$-gon ${\bf v}$, we define a non-degenerate pairing between the spaces ${\cal P}_{\bf v}$ and ${\cal P}_{\bf u}$ where ${\bf u} = {\bf v^\ast}$. Thus ${\cal P}_{\bf v^\ast} = ({\cal P}_{\bf v})^\ast$. The evolute map ${\cal P}_{\bf v} \to {\cal P}_{\bf v^\ast}$ is anti-self-adjoint.

The case of ``small-gons", that is, $(m+2)$- and $(m+3)$-gons in $\R^m$, is special. We prove that the second evolute of an  $(m+2)$-gons is homothetic to the original polygon (Theorem \ref{pentatheorem}) and that the first and third evolutes of an $(m+3)$-gon are homothetic (Theorem \ref{hexatheorem}). The first of these results is not new: it was proved by E. Tsukerman \cite{Tsu} by a different method. The second result was known in dimension 2: it was conjectured by B. Gr\"unbaum in \cite{Gr1,Gr2} and  proved in our paper \cite{ArFuIzTaTs}.

In general, if the  maximum module eigenvalue of the second evolute transformation is real, then the whole sequence of multiple evolutes of a polygon, considered up to parallel translations and scaling, is asymptotically 2-periodic (the limit coefficient of the homothety equals the maximum module eigenvalue; in particular, it may be negative). In Figures \ref{hepta8}, \ref{hepta9} and \ref{hepta1}, we show examples of such behavior for  heptagons.

Our last result on polygons concerns the spectrum of the second evolute map. Theorem \ref{spectrumdouble} states that every non-zero eigenvalue of this map has even multiplicity, generically, multiplicity 2. More precisely, the matrix of the restriction of the second evolute transformation to the image of the evolute transformation is conjugated to the block diagonal matrix of the form $\left[\displaystyle{\begin{array} {cc} A&0\\ 0&A^{\rm T}\end{array}}\right].$

The main ingredient of the proof is a construction of a linear symplectic structure on the space ${\cal P}_{\bf v}$, when $n-m$ is even, or on its quotient by the kernel of the evolute map, when $n-m$ is odd. The second evolute map is skew-Hamiltonian with respect to this symplectic structure.

Theorem \ref{spectrumdouble} is our strongest result on polygons; it provides alternative proofs of  Theorems \ref{pentatheorem} and \ref{hexatheorem}.

In the second part, we consider curves in $\R^3$ with non-vanishing curvature and non-vanishing torsion (the former is a general position property, whereas the latter is not). Our curves may have cusps, such as the curve $(t^2,t^3,t^4)$ at the origin, and they are equipped with the unit tangent vector field $T$, continuous through the cusps (this implies that the number of cusps is necessarily even). The vectors $T$ define the  tangent Gauss image of such a curve,  its tangent indicatrix, is an immersed closed locally convex curve on the unit sphere. 

The osculating sphere of a spacial curve passes through a quadruple of infinitesimally close points of the curve. The {\it evolute} of a curve is the locus of the centers of its osculating spheres (this evolute is also called the evolute of the 2nd kind in \cite{BL} and the focal curve in \cite{UrVa}). Equivalently, the evolute is the enveloping curve of the family of the normal planes of a spacial curve, and it is also called the edge of regression of the polar developable in \cite{St}.
See \cite{Fu,UrVa} for a study of evolutes of spacial curves. 

The space of curves in $\R^3$ modulo parallel translations fibers over the space of spherical curves, their tangent indicatrices. The space of curves ${\cal C}_{\gamma}$ with a fixed tangent indicatrix $\gamma\subset S^2$ is an infinite-dimensional vector space (for a precise definition of the space $C_\gamma$, see Section \ref{fixdtan}). One of our results is that the tangent indicatrix $\overline\gamma$ of the evolute of a spacial curve is spherically dual to $\gamma$, the indicatrix of the original curve, and the evolute map ${\cal C}_{\gamma} \to {\cal C}_{\overline\gamma}$ is linear (see Theorem \ref{linmap}). 

We prove continuous analogs of some of the results obtained for polygons. In particular, for a generic tangent indicatrix $\gamma$, we construct a non-degenerate pairing between the spaces ${\cal C}_{\gamma}$ and ${\cal C}_{\overline\gamma}$, and we prove that the evolute map ${\cal C}_{\gamma} \to {\cal C}_{\overline\gamma}$ is an anti-self-adjoint linear bijection.

The second evolute of a spacial curve has the same, up to a central symmetry, tangent indicatrix as the original curve, and one wonders whether there exist curves which are homothetic to their second evolutes. In the plane, the curves with this property are the classical hypocycloids, see \cite{ArFuIzTaTs}, Corollary 2.8.

 We construct a family of spacial curves that are homothetic to their second evolutes. The tangent indicatrices of these curves are circles on the unit sphere. These curves may be regarded as spacial analogs of hypocycloids, and indeed, they look somewhat like the classical hypocycloids, see Figure \ref{hypocycloids}. We do not know whether there are  other spacial curves homothetic to their second evolutes.

\paragraph{Acknowledgment.} We are grateful to A. Akopyan, A. Bobenko, I. Izmestiev, W. Schief, and E. Tsukerman for useful discussions, and especially to Yu. Suris for informing us about the results of \cite{FMMX}, that play the critical role in our proof of  Theorem \ref{spectrumdouble}.

The second author was  supported by NSF grant  DMS-1510055 and by the DFG Collaborative Research Center TRR 109 ``Discretization in Geometry and Dynamics".
The first and second authors gratefully acknowledge the hospitality of IHES and of Universit\"atszentrum Obergurgl, respectively.

We are grateful to the anonymous referee for the constructive criticism. 

\section{Polygons}\label{polygons} 
Let $n\geq m+2$. A (closed) $n$-{\it gon} in $\R^m$ is, by definition, a cyclically ordered sequence of lines with the following two properties:\smallskip

(1) every two consecutive lines have precisely one common point;

(2) no $m$ consecutive lines are parallel to (equivalently: are contained in) a hyperplane in $\R^m$.\smallskip

We call the lines and their intersection points the {\it sides} and the {\it vertices} of the polygon.

A common notation: the lines which comprise a polygon are given consecutive even or odd indices modulo $2n$: $\ell_2,\ell_4,\ldots,\ell_{2n}$ or $\ell_1,\ell_3,\ldots,\ell_{2n-1}$; the intersection point of $\ell_{i-1}$ and $\ell_{i+1}$ is denoted as $P_i$. A polygon is called {\it non-degenerate}, if no two consecutive vertices $P_i$ coincide. A non-degenerate polygon is determined by the sequence of its vertices, and, for such polygons, we can use the notation ${\bf P}=(P_1,P_3,\ldots, P_{2n-1})$ or $(P_2,P_4,\dots,P_{2n})$. In general, some consecutive vertices of a polygon can coincide (even all of them may be the same point). In such a degenerate case, we still may use the notation as above, having in mind that if $P_{i-1}=P_{i+1}$, then there must be specified a line $\ell_i$ passing through this point.

\subsection{Evolutes: definitions}\label{evoldef}  
There are two natural notions of the evolute of a polygon $\bf P$; following the terminology of \cite{ArFuIzTaTs}, we call them $\mathcal P$-{\it evolute} and $\mathcal A$-{\it evolute}. 

The $\mathcal P$-evolute $\bf Q$ of $\bf P$ is the $n$-gon whose (cyclically ordered) vertices are centers of the spheres passing through  $(m+1)$-tuples of (cyclically) consecutive vertices of $\bf P$. In the degenerate case, the condition of passing through coinciding vertices $P_{j-1}, P_{j+1}$ should be supplemented by the condition of tangency to the side $\ell_i$. The most convenient way of numeration of the vertices of the evolutes is as follows. The vertex $Q_i$ of $\bf Q$ is  the center of the sphere passing through the points$$P_{i-m},P_{i-m+2},P_{i-m+4},\ldots,P_{i+m-4},P_{i+m-2},P_{i+m}$$(we consider the subscripts as defined modulo $2n$). [Thus, if $m$ is even (odd), then the vertices of $\bf Q$ are numerated as vertices (sides) of $\bf P$.] The point $Q_i$ may be also described as the point of intersection of the perpendicular bisector hyperplanes of the sides $\ell_{i-m+1},\ell_{i-m+3},\ldots,\ell_{i+m-3},\ell_{i+m-1}$ of $\bf P$ (in the case $P_{j-1}=P_{j+1}$, the perpendicular bisector hyperplane to the side $\ell_j$ is the plane through the point $P_{j-1}=P_{j+1}$ perpendicular to $\ell_j$). The last description of $Q_i$ shows that the side $q_j=\overline{Q_{j-1}Q_{j+1}}$ of $\bf Q$ is perpendicular to the sides $\ell_{j-m+2},\ell_{j-m+4},\ldots,\ell_{j+m-4},\ell_{j+m-2}$; in the degenerate case $Q_{j-1}=Q_{j+1}$ this property of the side $q_j$ should be taken for its definition.

We will define the $\mathcal A$-evolute $\bf R$ of $\bf P$ only in the case when both $\bf P$ and $\bf R$ (as defined below) are non-degenerate. $\bf R$ is the $n$-gon whose vertices are the centers of the spheres  tangent to the $(m+1)$-tuples of (cyclically) consecutive sides of $\bf P$. The vertices of $\bf R$ are numerated as follows. The vertex $R_i$ of $\bf R$ is the center of the sphere tangent to the sides $\ell_{i-m},\ell_{i-m+2},\dots,\ell_{i+m-2},\ell_{i-m}$ of $\bf P$. [Thus, if $m$ is even (odd), then the vertices of $\bf R$ are numerated as sides (vertices) of $\bf P$.] The point $R_i$ may be described as the point of intersection of the bisectorial planes of the angles $$\widehat{P_{i-m-1}P_{i-m+1}P_{i-m+3}},\ \widehat{P_{i-m+1}P_{i-m+3}P_{i-m+5}},\ldots,\widehat{P_{i+m-3}P_{i+m-1}P_{i+m+1}}$$ of $\bf P$. Since every angle has two different bissectorial planes, an $n$-gon has, in general, a large amount of different $\mathcal A$-evolutes. 

\begin{proposition}\label{secondev} The sides of the second $\mathcal P$-evolute of a polygon $\bf P$ are parallel to the sides of the polygon $\bf P$.\end{proposition}

\proof Let $s_i$ be sides of the $\mathcal P$-evolute $\bf S$ of the $\mathcal P$-evolute $\bf Q$ of $\bf P$. The sides $q_{i-m+2},q_{i-m+4},\ldots,q_{i+m-4},q_{i+m-2}$ of $\bf Q$ are contained in the perpendicular bisector hyperplane of the side $\ell_i$ of $\bf P$. Likewise, the side $s_i$ is contained in the perpendicular bisector planes of the sides $q_{i-m+2},q_{i-m+4},\ldots,q_{i+m-4},q_{i+m-2}$. It follows that the side $\ell_i$ of $\bf P$ and the side $s_i$ of $\bf S$ are both perpendicular to the sides $q_{i-m+2},q_{i-m+4},\ldots,q_{i+m-4},q_{i+m-2}$ of $\bf Q$; hence they are parallel to each other. (The $m-1$ sides listed above are not parallel to an $(m-2)$-dimensional plane.) \proofend

\subsection{Involutes: existence and uniqueness}\label{involutes} 
If an $n$-gon $\bf Q$ is the $\mathcal P$-evolute of an $n$-gon $\bf P$, then we call $\bf P$ a $\mathcal P$-involute of $\bf Q$. If an $n$-gon $\bf Q$ is one of $\mathcal A$-evolutes of an $n$-gon $\bf P$, then we call $\bf P$ an $\mathcal A$-involute of $\bf Q$.

\begin{lemma}\label{fixedpts} 
The table below contains a complete information on the  existence and uniqueness of fixed points and invariant lines for a generic\footnote{The word {\sl generic} means that the statement holds within a dense open set of isometries.} isometry $\sigma\colon\R^m\to\R^m$:

\begin{figure}[hbtp] 
\centering
\includegraphics[width=5.4in]{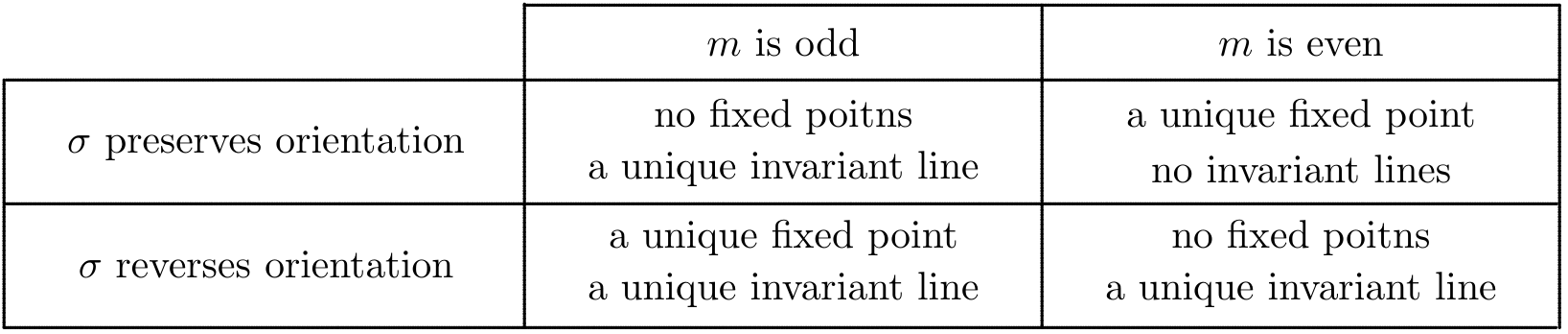}
\label{tablica}
\end{figure}

\end{lemma}
\vskip.2in

\proof The isometry $\sigma$ is a composition of an orthogonal transformation $\sigma_0$ and a parallel translation. 

Let $m$ be even and $\sigma$ be orientation preserving. Generically, $\pm1$ is not an eigenvalue of $\sigma_0$. This means that $\sigma$ has no invariant lines (even no invariant directions). 
Let us prove, by induction, that if 1 is not an eigenvalue of $\sigma_0$, then $\sigma$ has a unique fixed point (for $m=2$ the fact is obvious). 

Let $M_0$ be an invariant 2-dimensional plane of $\sigma_0$ such that the action of $\sigma_0$ in $M_0$ is orientation preserving (this must exist). Then every 2-dimensional plane parallel to $M_0$ is mapped by $\sigma$ into a plane parallel to $M_0$, and $\sigma$ is factorized to an isometry $\widetilde\sigma$ of $\R^m/M_0$. The orthogonal transformation $\widetilde\sigma_0$ also has no eigenvalues 1, so, by the induction hypothesis, $\widetilde\sigma$ has a unique fixed point. Hence, there exists a unique $\sigma$-invariant plane $M\subset\R^m$ parallel to $M_0$, and every fixed point of $\sigma$ must be contained in $M$. The restriction of $\sigma$ to $M$ is a rotation by a non-zero angle, it has a unique fixed point, which is a unique fixed point of $\sigma$.

If $m$ is odd and $\sigma$ preserves orientation, then 1 is an eigenvalue of $\sigma_0$, generically, of multiplicity 1. Hence, there exists a line $L_0$ which is pointwise fixed by $\sigma_0$. By the previous case, the arising isometry $\widetilde\sigma$ of $\R^m/L_0$ has a unique fixed point, so there exists a unique $\sigma$-invariant line $L\subset\R^m$ parallel to $L_0$. If $\sigma$ has a fixed point, then this point must be contained in $L$, but the restriction of $\sigma$ to $L$ is a parallel translation, and it is generically fixed point free.

If $m$ is odd and $\sigma$ reverses orientation, then $\sigma_0$ has an eigenvalue $-1$ and, generically, no eigenvalue 1. Hence, there is a $\sigma_0$-invariant line $L_0$, on which $\sigma_0$ acts as a reflection in 0. The arising isometry $\sigma_0$ of $\R^m/L_0$ preserves  orientation and 1 is not an eigenvalue of $\widetilde\sigma_0$; hence, $\widetilde\sigma$ has a unique fixed point. Thus there is a unique $\sigma$-invariant line $L\subset\R^m$, and the direction of $L$ is reversed by $\sigma$. Hence, $L$ contains a unique fixed point of $\sigma$.

Finally, if $m$ is even and $\sigma$ reverses orientation, then both 1 and $-1$ are eigenvalues of $\sigma_0$, generically, both of multiplicity 1. Then there is a $\sigma_0$-invariant plane $M_0\subset\R^m$, and the restriction of $\sigma_0$ to this plane is a reflection in a line. The arising isometry $\widetilde\sigma$ of $\R^m/M_0$ preserves  orientation, so it has a unique fixed point. Thus, there is a unique $\sigma$-invariant plane $M\subset\R^m$, and the restriction of $\sigma$ to this plane reverses  orientation. Generically, it is a glide reflection. Hence, $M$ contains a unique $\sigma$-invariant line and no fixed points. \proofend

Using the lemma, we address the question of the existence and uniqueness of involutes.

\begin{proposition}\label{existinvol} If $n-m$ is odd, then a generic\footnote{in particular, non-degenerate} $n$-gon has no $\mathcal P$-involutes. If $n-m$ is even, then a generic $n$-gon has a unique $\mathcal P$-involute. For any $n\ge m+2$, a generic $n$-gon has a unique $\mathcal A$-involute.\end{proposition}

\proof In this proof we assume that $m$ is odd; the proof for even $m$ is the same, up to some modification of notations (for example, the points $P_i$ will be numerated by even numbers, while the reflections $s_i$ will be numerated by odd numbers). 

Let ${\bf Q}=(Q_2,Q_4,\ldots Q_{2n})$ be an $n$-gon satisfying the conditions formulated at the beginning of this section. Let $s_i, i=2,4,\dots,2n$, be the reflection of $\R^m$ in the $(m-1)$-dimensional  plane containing the $m$ points $Q_{i-m+1},Q_{i-m+3},\dots,Q_{i+m-3},Q_{i+m-1}$, and let $\sigma=s_{2n}\circ\dots\circ s_4\circ s_2$. Obviously, for a generic $n$-gon, $\sigma$ is a generic isometry, preserving or reversing orientation, depending on the parity of $n$. 

The polygon $\bf Q$ is a $\mathcal P$-evolute of a polygon ${\bf P}=(P_1,P_3,\dots,P_{2n-1})$ if and only if, for every $i=2,4,\dots,2n$, the point $Q_i$ belongs to the perpendicular bisector hyperplanes of the sides $\ell_{i-m+2},\ell_{i-m+4},\ldots,\ell_{i+m-4},\ell_{i+m-2}$ of $\bf P$. In other words, for every $i$, the points $Q_{i-m+1},Q_{i-m+3},\dots,Q_{i+m-3},Q_{i+m-1}$ lie in the perpendicular bisector hyperplane of $\ell_i$, which means, in turn, that $P_{i+1}=s_i(P_{i-1})$. Thus, $\sigma(P_1)=P_{2n+1}=P_1$. 

If $n$ is even, then $\sigma$ preserves orientation, and by Lemma \ref{fixedpts}, it has no fixed points, so $\bf P$ cannot exist. If $n$ is odd, then $\sigma$ has a generically unique fixed point. We take this point for $P_1$ and put $P_3=s_2(P_1),\, P_5=s_4(P_3),\, P_7=s_6(P_5)$, and so on; we obtain an $n$-gon $\bf P$, whose $\mathcal P$-evolute is $\bf Q$. 

The proof for $\mathcal A$-involutes is similar, the only difference is that we begin not with a fixed point, but with an invariant line of $\sigma$. We successively apply to this line the reflections $s_2,s_4,s_6,\ldots$. We obtain $n$ lines which form an $n$-gon  whose $\mathcal A$-evolute is $\bf Q$. (Notice that if we reflect a line in a hyperplane not parallel to this line, then the line and its reflection intersect at the intersection point of the line with the hyperplane.) \proofend

In the rest of the paper, we will not consider either $\mathcal A$-evolutes or $\mathcal A$-involutes. So, we will refer to $\mathcal P$-evolutes and $\mathcal P$-involutes simply as to evolutes and involutes.

\subsection{Spherical polygons and polygons with fixed directions of the sides}\label{TvPv}

A spherical $n$-gon is a cyclically ordered collection of $n$ points of the unit sphere $S^{m-1}$; notation: ${\bf v} = (v_2,v_4,\ldots,v_{2n})$ or $(v_1, v_3,\ldots, v_{2n-1})$. We assume that the spherical polygons are generic: no (cyclically) consecutive $m$-tuple of vectors $v_i$ is linearly dependent.

The {\it signature} of a polygon $\bf v$ is the cyclic sequence $\bf s$ of signs $s_i$ of the determinants $d_i=\det(v_{i-m+1},v_{i-m+3},\dots,v_{i+m-3},v_{i+m-1})$.

For a spherical $n$-gon ${\bf v} = (v_2,v_4,\ldots,v_{2n})$, let ${\bf u} = {\bf v}^\ast$ be the {\it dual} $n$-gon ${\bf u} = (u_1,u_3,\ldots,u_{2n-1})$ (${\bf u}=(u_2,u_4,\dots,u_{2n})$, if $m$ is even) where $u_i$ is the positive unit normal vector of the  $(m-1)$-dimensional plane spanned by the vectors $v_{i-m+2},v_{i-m+4},\dots,v_{i+m-4},v_{i+m-2}$. The positivity means that the sign of $u_i\cdot v_{i+m}$ is $s_{i+1}$, and the sign of $u_i\cdot v_{i-m}$ is $(-1)^{m-1}s_{i-1}$; in particular, both these dot products are not zero.

\begin{lemma}\label{dualgeneric} If $\bf v$ is generic, then ${\bf u}={\bf v}^\ast$ is generic.\end{lemma}

\proof Suppose that \begin{equation}\label{lindep}A_{i-m+1}u_{i-m+1}+A_{i-m+3}u_{i-m+3}+\dots+A_{i+m-1}u_{i+m-1}=0.\end{equation} Dot (\ref{lindep}) with $v_{i+1}$; we get $A_{i-m+1}u_{i-m+1}\cdot v_{i+1}=0$, and since $u_{i-m+1}\cdot v_{i+1}\neq0$, we have $A_{i-m+1}=0$. Taking this into account, dot (\ref{lindep}) with $v_{i+3}$; we get $A_{i-m+3}u_{i-m+3}\cdot v_{i+3}=0$, and since $u_{i-m+3}\cdot v_{i+3}\neq0$, we have  $A_{i-m+3}=0$. And so on. \proofend

Let ${\bf w}=(w_2,w_4,\dots,w_{2n})$ be $({\bf v}^\ast)^\ast$.

\begin{lemma}\label{plusminus} One has $w_i=\pm v_i$.\end{lemma}

\proof Indeed, both $v_i$ and $w_i$ are orthogonal to $m-1$ linearly independent vectors $v_{i-m+2},v_{i-m+4},\dots,v_{i+m-4},v_{i+m-2}$. \proofend

The following two lemmas provide clarification to Lemmas \ref{dualgeneric} and \ref{plusminus}. We will not use these statements, and we leave their proofs, which may be regarded as exercises in linear algebra, to the reader.

\begin{lemma}\label{vstarstar} One has $w_i=(-1)^{m-1}s_{i-m+3}s_{i-m+5}\dots s_{i+m-5}s_{i+m-3}v_i.$ \proofend\end{lemma}

Let $\{s^\ast_i\}$ be the signature of the dual polygon $\bf u$.

\begin{lemma}\label{vstarsign} One has $s^\ast_i=s_{i-m+2}s_{i-m+4}\dots s_{i+m-4}s_{i+m-2}.$ \proofend\end{lemma}
 
Fix a spherical $n$-gon $\bf v$, and let ${\mathcal P}_{\bf v}$ be the space of $n$-gons in $\R^m$, considered up to a parallel translation, whose sides are parallel to the respective vectors $v_i$. If ${\bf P} = (P_1, P_3,\ldots, P_{2n-1})$ is such a polygon, define the real numbers $x_i$ (signed side lengths) by 
$$P_{i+1}-P_{i-1}=x_iv_i.$$
The vector ${\bf x} = (x_2,\ldots,x_{2n})$ uniquely determines the polygon (up to parallel translation). 

Notice that some (or even all) signed side lengths $x_i$ may be zero. In other words, we do not exclude the possibility that the polygon is degenerate, but when we consider a polygon as belonging to ${\mathcal P}_{\bf v}$, we assign directions to all the sides.

The coordinates $x_i$ satisfy  $m$ linear relations
$
\sum x_i v_i =0,
$
saying that the respective polygon closes up. Thus ${\mathcal P}_{\bf v}$ is a vector space of dimension $n-m$.
Notice that the space  ${\mathcal P}_{\bf v}$ remains the same if we replace some vectors $v_i$ by the opposite vectors (although some coordinates $x_i$ change signs).

\subsection{The evolute transformation}\label{evtransfpoly}

In this section we show that the evolute of a polygon ${\bf P}\in{\mathcal P}_{\bf v}$ belongs to ${\mathcal P}_{\bf u}$, where $\bf u=v^\ast$, and that the evolute map ${\mathcal E}\colon{\mathcal  P}_{\bf v}\to{\mathcal P}_{\bf u}$ is a linear transformation.

\begin{lemma} \label{pevol}
Let ${\bf P}\in{\mathcal P}_{\bf v}$ and ${\bf Q} = {\cal E} ({\bf P}).$ Then
 ${\bf Q} \in {\cal P}_{\bf u}$, where ${\bf u} = {\bf v}^*$.
\end{lemma}

\proof
By continuity, it is sufficient to prove the statement in the case when both $\bf P$ and $\bf Q$ are non-degenerate. The point $Q_{i-1}$ has equal distances to the points$$P_{i-m-1},P_{i-m+1},\dots,P_{i+m-3},P_{i+m-1},$$ and the point $Q_{i+1}$ has equal distances to the points$$P_{i-m+1},P_{i-m+3},\dots,P_{i+m-1},P_{i+m+1}.$$ Hence, both these points have equal distances to the points$$P_{i-m+1},P_{i-m+3},\dots,P_{i+m-3},P_{i+m-1},$$ and the vector $Q_{i+1}-Q_{i-1}$ is orthogonal to the vectors $P_{i-m+3}-P_{i-m+1},$ $\dots,P_{i+m-1}-P_{i+m-3}$. This means  that if ${\bf P}\in{\mathcal P}_{\bf v}$, then ${\bf Q}\in{\mathcal P}_{{\bf v}^\ast}$. \proofend

Let $a_1,\dots,a_{m+1}$ be $m+1$ points in $\R^m$, and let $z$ be their circumcenter.

\begin{lemma}\label{circumformula} 
The point $z$ is determined by the system of linear equations$$z\cdot(a_{i+1}-a_1)=\frac{a_{i+1}\cdot a_{i+1}-a_i\cdot a_i}2,\, i=1,\ldots,m.$$
\end{lemma}

\proof The circumcenter $z$ satisfies the equations$$|z-a_1|^2 =\ldots=|z-a_{m+1}|^2,$$that are equivalent to the above system of linear equations. \proofend

Let $\bf P\in{\mathcal P}_v$, and let $\bf Q={\mathcal E}(P)\in{\mathcal P}_u$. The vertices $P_j$ of $\bf P$ and $Q_i$ of $\bf Q$ are labeled by odd or even integers depending on the parity of $m$. Define the real numbers (signed side lengths) $y_{i}$ by 
$$Q_{i+1}-Q_{i-1} = y_iu_i.$$
We obtain a vector ${\bf y} = \{y_i\}$, and the evolute transformation ${\mathcal E}\colon {\mathcal P}_{\bf v}\to{\mathcal P}_{\bf u}$ acts by $\bf x\mapsto y$.

\begin{theorem}\label{linear} 
The vector $\bf y$ is obtained from $\bf x$ by a linear transformation depending on $\bf v$.
\end{theorem}

\proof Again, it is sufficient to prove the linearity for a non-degenerate $\bf P$. Consider $m+2$ consecutive vertices $P_{i-m-1},P_{i-m+1},\ldots,P_{i+m-1},$ $P_{i+m+1}$
of {\bf P}; the $m+1$ consecutive vectors of the sides are$$x_{i-m}v_{i-m},\, x_{i-m+2}v_{i-m+2},\ldots, x_{i+m-2}v_{i+m-2},x_{i+m}v_{i+m}.$$Assume that the first vertex is the origin. Then the $m+2$ vertices are
\begin{equation} \label{eq8}
P_{i-m+2j-1}=\sum_{k=0}^{j-1} x_{i-m+2k}v_{i-m+2k},\ j=0,\ldots,m+1.
\end{equation}
The two relevant vertices of the evolute are $Q_{i-1}$ and $Q_{i+1}$, and the relevant unit vector is $u_i.$ Recall that $Q_{i-1}$ is the circumcenter of the $m+1$ points $P_{i-m+2j-1},\, j=0,\ldots,m$, and  $Q_{i+1}$ is the circumcenter of the $m+1$ points $P_{i-m+2j-1},\, j=1,\ldots,m+1$.

Apply Lemma \ref{circumformula} to $Q_{i-1}\cdot\left(P_{i-m+2j+1}-P_{i-m+2j-1}\right),\, j=0,\ldots,m$, to obtain
$$Q_{i-1}\cdot\left(x_{i-m+2j}v_{i-m+2j}\right)=\frac{|P_{i-m+2j+1}|^2-|P_{i-m+2j-1}|^2}2,$$which shows, in virtue of (\ref{eq8}), that$$Q_{i-1}\cdot v_{i-m+2j}=\frac{2\left(x_{i-m}v_{i-m}+\ldots+x_{i-m+2j-2}v_{i-m+2j-1}\right)\cdot v_{i-m+2j}+x_{i-m+2j}}2.$$We see that $Q_{i-1}$ is a linear function of $x_{i-m+2j},\, j=0,\dots,m$; similarly, $Q_{i+1}$ is a linear function of $x_{i-m+2j},\, j=1,\dots,m+1$. Thus, $Q_{i+1}-Q_{i-1}$ is a linear function of $x_{i-m+2j},\, j=0,\dots,m+1$, and so is $y_i=(Q_{i+1}-Q_{i-1})\cdot u_i$. \proofend

\subsection{Rank of the evolute transformation}\label{rank} 

The following proposition describes the rank of the evolute transformation.

\begin{proposition}\label{kernel} 
Generically, the evolute transformation ${\mathcal E}\colon{\mathcal P}_{\bf v}\to{\mathcal P}_{\bf u}$ has 1-dimensional kernel if $n-m$ is odd, and has full rank if $n-m$ is even.
\end{proposition}

\proof 
 If $n-m$ is even, then the map ${\mathcal E}\colon{\mathcal P}_{\bf v}\to{\mathcal P}_{\bf u}$ has a full rank, because it is onto: in this case, a generic $n$-gon has an involute (Proposition \ref{existinvol}). 
 
Let $n-m$ be odd. We will show that $\dim\Ker{\mathcal E}\geq1$ and, generically, this dimension is 1. A polygon belongs to the kernel of the evolute map if and only if it is inscribed into a sphere. We will show that, for every $\bf v$, the space ${\mathcal P}_{\bf v}$ contains a polygon inscribed into the unit sphere $S^{m-1}$ and, generically, this polygon is unique up to the reflection in the center. 
 
Take a point $P_1\in S^{m-1}$ and reflect it, successively, in the hyperplanes passing through the origin and perpendicular to the vectors $v_2,v_4,\ldots,v_{2n}$; we obtain points $P_3, P_5,\ldots, P_{2n+1}\in S^{m-1}$. If $P_{2n+1}=P_1$, then ${\bf P}=(P_1,P_3,\ldots,P_{2n-1})\in {\mathcal P}_{\bf v}$ is  inscribed in $S^{m-1}$, and all polygons from ${\mathcal P}_{\bf v}$ inscribed in $S^{m-1}$ are obtained by this construction. 

The transformation $P_1\mapsto P_{2n+1}$ is an isometry of $S^{m-1}$ of degree $(-1)^n$. The Lefschetz number of this map is 2, so it has fixed points. Generically, there are precisely 2 fixed points, and they are antipodal. Thus, the space of inscribed $n$-gons is not zero, and generically has dimension 1. \proofend

\subsection{${\mathcal P}_{\bf v}$ as the dual space to ${\mathcal P}_{\bf v^\ast}$}\label{dualsect}

Consider an $n$-gon $\bf P\in{\mathcal P}_v$, where ${\bf v}=\{v_i\}\subset S^{m-1}$ is a generic spherical $n$-gon. Choose an origin and  define the support number $\lambda_i$  as the signed distance from the origin to the $(m-1)$-dimensional plane which contains the sides $\ell_{i-m+2},\ell_{i-m+4},\ldots,\ell_{i+m-4},\ell_{i-m+2},$ parallel to the vectors $v_{i-m+2},v_{i-m+4},\ldots,v_{i+m-4},v_{i-m+2}$ and oriented by these vectors. (The polygon $\bf P$ is not supposed to be generic; for example, some or all its vertices may coincide.)

Let $\bf Q \in {\mathcal P}_{\bf v^\ast}$. Define $\langle{\bf P,Q}\rangle=\sum_iy_i\lambda_i$ where $y_i$ is the signed length of the $i$-th side of $\bf Q$.

\begin{lemma} \label{welldef}
The pairing $\langle{\bf P,Q}\rangle$ is well defined, that is, it does not depend on the choice of the origin. 
\end{lemma}

\proof
The orienting normal vector of the hyperplane parallel to the vectors $v_{i-m+2}, v_{i-m+4},\ldots,v_{i+m-4},v_{i-m+2}$ is $u_i$, therefore $\lambda_i=P_{i+k}\cdot u_i$ for all $k= -m+1, -m+3,\ldots, m-3,m-1$. It follows that $y_i\lambda_i=P_{i+k}\cdot(Q_{i+1}-Q_{i-1})$ for the same values of $k$, and hence 
$$
\langle{\bf P,Q}\rangle = \sum_i P_{i+k}\cdot(Q_{i+1}-Q_{i-1}).
$$
A parallel translation through vector $R$ changes this expression by
$$
\sum_i R\cdot(Q_{i+1}-Q_{i-1}) = R\cdot\sum_i \left(Q_{i+1}-Q_{i-1}\right)=0, 
$$
as claimed.
\proofend

The pairing $\langle\ \rangle\colon{\mathcal P}_{\bf v}\otimes{\mathcal P}_{\bf v^\ast}\to\R$ is skew-symmetric in the following sense.

\begin{proposition} \label{antisym} For any ${\bf P}\in{\mathcal P}_{\bf v},\ {\bf Q}\in{\mathcal P}_{\bf v^\ast}$, one has $\langle{\bf P,Q}\rangle=-\langle{\bf Q,P}\rangle$. 
\end{proposition}

\proof If $m$ is odd, then 
$$
\langle{\bf P,Q}\rangle=\sum_i P_i\cdot(Q_{i+1}-Q_{i-1})=\sum_i Q_i\cdot(P_{i-1}-P_{i+1})=-\langle{\bf Q,P}\rangle,$$ and if $m$ is even, then 
$$
\langle{\bf P,Q}\rangle=\sum_i P_{i+1}\cdot(Q_{i+1}-Q_{i-1})=\sum_i Q_{i-1}\cdot(P_{i-1}-P_{i+1})=-\langle{\bf Q,P}\rangle,$$
as claimed. \proofend

\begin{proposition} \label{nondeg} The pairing $\langle\ \rangle\colon{\mathcal P}_{\bf v}\otimes{\mathcal P}_{\bf v^\ast}\to\R$ is non-degenerate. Thus, with respect to this pairing, ${\mathcal P}_{\bf v^\ast}$ is  the dual space $\left({\mathcal P}_{\bf v}\right)^\ast$.\end{proposition}

\proof We want to prove that for every non-zero set $\{\lambda_i\}$, there exists a set $\{y_i\}$, satisfying the condition $\sum_iy_iu_i=0$, such that $\sum_iy_i\lambda_i\ne0$. 

The opposite would mean that $\sum_iy_iu_i=0$ implies $\sum_iy_i\lambda_i=0$, that is, $\lambda_i$ is a linear combination of the coordinates of $u_i$ with coefficients not depending on $i$, $\lambda_i=C\cdot u_i$ for some constant vector $C$. But $\lambda_i=P_{i+k}\cdot u_i$, for $k=-m+1,-m+3,\ldots,m-3,m-1$. Thus $(P_{i+k}-C)\cdot u_i=0$ for the same values of $k$. 

By applying a parallel translation, we may assume that $C=0$. 
The orthogonal complement of the vector $u_i$ is spanned by the vectors $P_{i-m+3} - P_{i-m+1},  P_{i-m+5} - P_{i-m+3}, \ldots, P_{i+m-1} - P_{i+m-3}$, implying a linear relation between the $m$ vectors $P_{i+k}, k=-m+1,-m+3,\ldots,m-3,m-1$. This holds for every $i$, hence the polygon $\bf P$ lies in a hyperplane, contradicting  the genericity of $\bf v$.  \proofend

The evolute transformation ${\mathcal E}:{\mathcal P}_{\bf v}\to{\mathcal P}_{\bf v^\ast}=\left({\mathcal P}_{\bf v}\right)^\ast$ is anti-self-adjoint. 

\begin{proposition} \label{orthogonal} For every ${\bf P}\in{\mathcal P}_{v^\ast}$,  one has  $\langle{\bf P,{\mathcal E}(P)}\rangle=0$.\end{proposition}

\proof Let $\bf Q={\mathcal E}(P)$. If $m$ is odd, then $\langle{\bf Q,P}\rangle=\sum_i Q_i\cdot(P_{i+1}-P_{i-1})$, and according to Lemma \ref{circumformula}, the last expression equals  $\sum\limits_i\displaystyle\frac{|P_{i+1}|^2-|P_{i-1}|^2}2=0$. For $m$ even, the proof is the same, only instead of $Q_i$ one takes $Q_{i+1}$. \proofend

\subsection{The cases of ``small-gons": first, second and third evolutes of ($\bf{m+2}$)- and ($\bf{m+3}$)-gons in $\bf\R^m$}\label{smallgons}

Since the assumption $n\geq m+2$ was made in the very beginning of Section \ref{polygons}, we start with the case $n=m+2$.

The following statement was proved by E. Tsukerman \cite{Tsu}. 

\begin{theorem} \label{pentatheorem} 
For a generic $(m+2)$-gon $\bf P$ in $\R^m$, its second evolute ${\mathcal E}^2({\bf P})$ is homothetic to $\bf P$.
\end{theorem}

(Actually, Tsukerman's work contains similar results for polygons in all spaces of constant curvature; we restrict ourselves here to the Euclidean case.) We will show that Theorem \ref{pentatheorem} is an almost immediate corollary of Proposition \ref{orthogonal}.

\paragraph{Proof\hskip6pt of\hskip6pt  Theorem\hskip6pt  \ref{pentatheorem}.} Let $\bf P$ be a non-zero element of some ${\mathcal P}_{\bf v}$. Notice that $\dim{\mathcal P}_{\bf v}=2$. By Proposition \ref{kernel}, the transformation ${\mathcal E}\colon{\mathcal P}_{\bf v}\to{\mathcal P}_{\bf v^\ast}$ has full rank, so ${\mathcal E}({\bf P})\ne0$. By Propositions \ref{antisym}, \ref{nondeg}, and \ref{orthogonal}, both $\bf P$ and $\bf{\mathcal E}^2(P)$ belong to the orthogonal complement of $\bf{\mathcal E}(P)\subset{\mathcal P}_{v^\ast}$, and this orthogonal complement is 1-dimensional. Thus, $\bf P$ and $\bf{\mathcal E}^2(P)$ are proportional, that is, the polygons $\bf P$ and $\bf{\mathcal E}^2(P)$ are homothetic. \proofend\smallskip

Now, let us turn to $(m+3)$-gons.

\begin{theorem} \label{hexatheorem} 
For a generic $(m+3)$-gon $\bf P$ in $\R^m$, its third evolute ${\mathcal E}^3({\bf P})$ is homothetic to its first evolute ${\mathcal E}(\bf P)$.
\end{theorem}

(For $m=2$, this is Proposition 4.15 of \cite{ArFuIzTaTs}.)

\proof Let $\bf P$ be a non-zero element of some ${\mathcal P}_{\bf v}$. By Proposition \ref{kernel}, the evolute transformation ${\mathcal E}\colon{\mathcal P}_{\bf v}\to{\mathcal P}_{\bf v^\ast}$ has a 1-dimensional kernel $\mathcal K$. The orthogonal complement ${\mathcal K}^\perp\subset{\mathcal P}_{\bf v^\ast}$ is the image of the evolute transformation. We may assume that ${\mathcal E}^2({\bf P})\notin{\mathcal K}$ (otherwise, ${\mathcal E}^3({\bf P})=0$, and the statement of the theorem is trivial). Thus, both ${\mathcal E}({\bf P})$ and ${\mathcal E}^3({\bf P})$ belong to two different two-dimensional subspaces ${\mathcal K}^\perp$ and ${\mathcal E}^2({\bf P})^\perp$ of the three-dimensional space ${\mathcal P}_{\bf v^\ast}$. The intersection of these subspaces is one-dimensional, hence the polygons ${\mathcal E}({\bf P})$ and ${\mathcal E}^3({\bf P})$ are homothetic. \proofend

\subsection{Iteration of the evolute transformations on general polygons} \label{iterations}

The second evolute transformation ${\mathcal P}_{\bf v}\to{\mathcal P}_{\bf v}$ is linear, and it is a standard fact from linear algebra that the asymptotic behavior of the sequence of iterations of this map depends largely on the maximal module eigenvalue of this transfiormation. If this eigenvalue is real, then sufficiently distant members of this sequence will be almost homothetic to each other. 

\begin{figure}[hbtp] 
\centering
\includegraphics[width=4.4in]{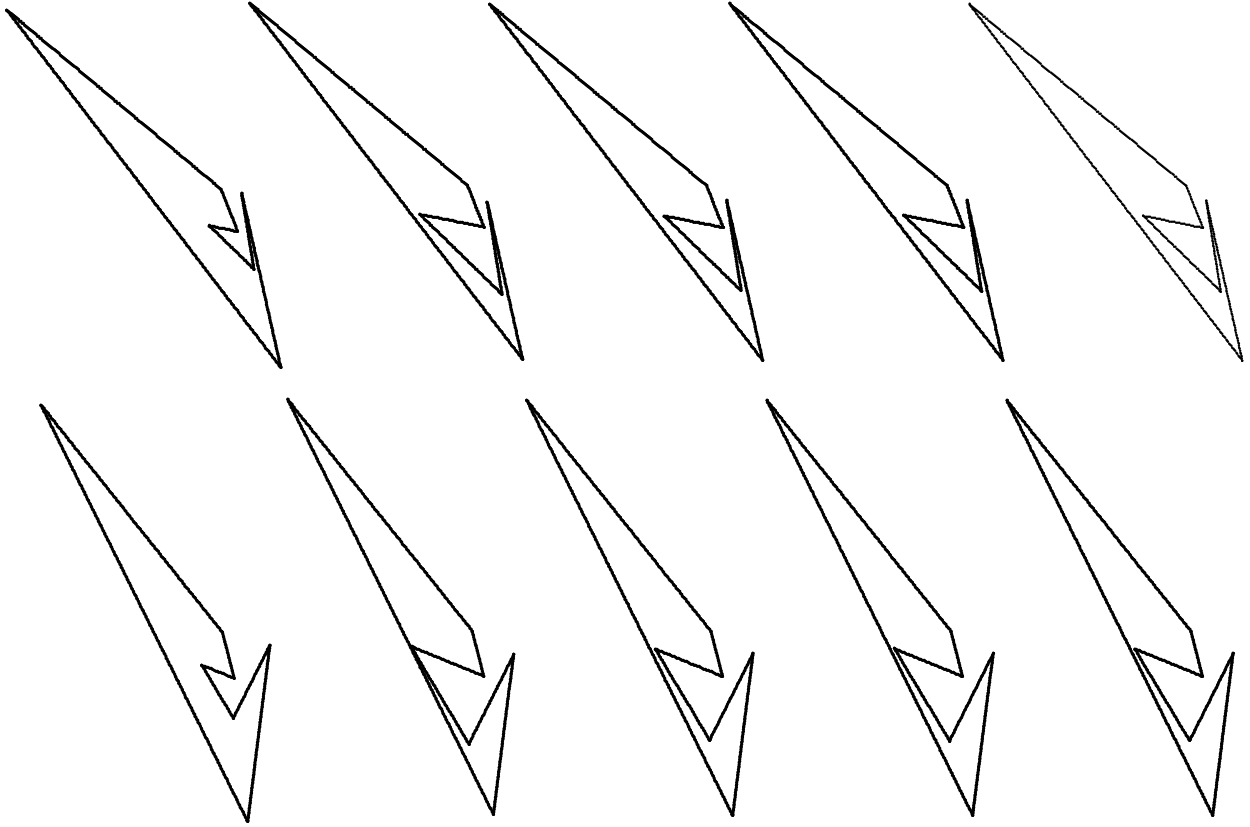}
\caption{Even-numbered evolutes of a heptagon in $\R^3$, two projections.}
\label{hepta8}
\end{figure}

To make the described property of evolutes better visible, we apply to each member of the sequence two transformations; first, a parallel translation which carries the centroid of the vertices to the origin, and a homothety centered at the origin, which makes the maximal distance from the origin to the vertices equal to 1. 

The appearance of the resulting sequence of polygons depends on the maximal eigenvalue of the transformation $\mathcal E^2$. 

\begin{figure}[hbtp] 
\centering
\includegraphics[width=4.4in]{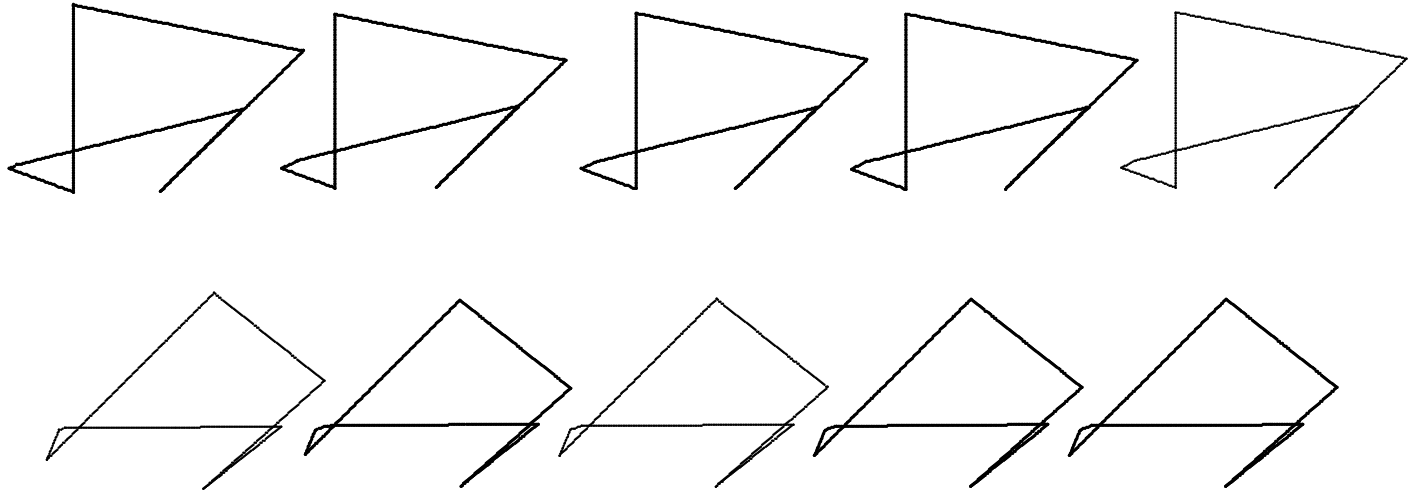}
\caption{Odd-numbered evolutes of the same heptagon in $\R^3$, two projections.}
\label{hepta9}
\end{figure}

If this eigenvalue is real and positive, then, starting with some moment, the polygons will be almost undistinguishable from each other. This case is presented in Figures \ref{hepta8} and \ref{hepta9}.

\begin{figure}[hbtp] 
\centering
\includegraphics[width=4.4in]{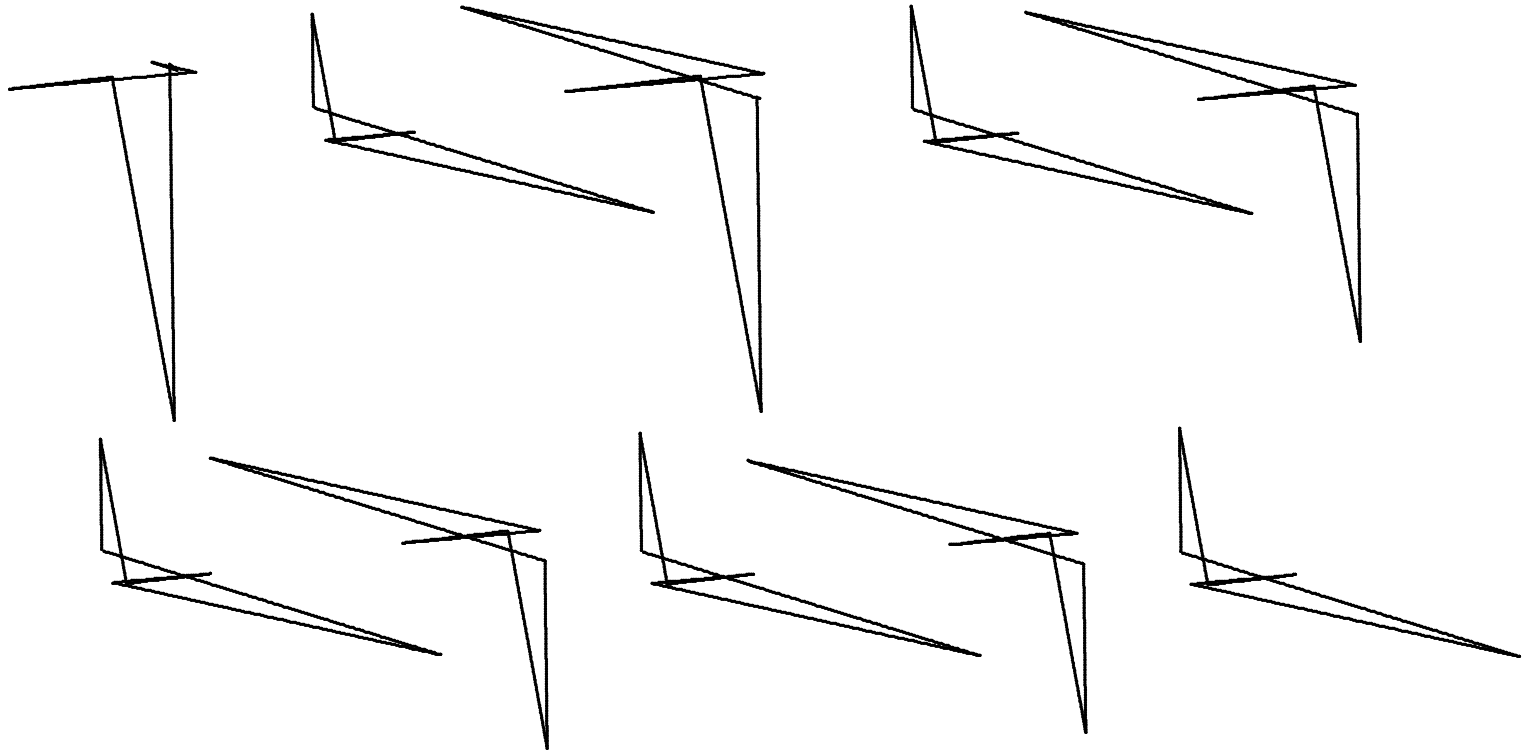}
\caption{The sequence of even-numbered evolutes of a heptagon in the case when the maximum module eigenvalue of the transformation $\mathcal E^2$ is negative (one pojection).}
\label{hepta1}
\end{figure}

If the maximum module eigenvalue is real and negative, then the shapes of the distant members of the sequence will be almost the same, but the transition from a polygon to the next one will include a flip: the reflection in the origin (this case is presented in Figure \ref{hepta1}). 

Finally, if the maximum module eigenvalue is not real, then the sequence of evolutes may not reveal periodicity (generically, one observes a quasi-periodic behavior).

The second evolute transformation acts separately on even numbered and odd numbered evolutes of a polygon. Let ${\bf v}$ be a spherical polygon and ${\bf u}$ its dual. 

\begin{proposition} \label{conj}
The second evolute maps ${\mathcal P}_{\bf v}\to{\mathcal P}_{\bf v}$ and ${\mathcal P}_{\bf u}\to{\mathcal P}_{\bf u}$ have the same eigenvalues.
\end{proposition}

\proof One has two first evolute maps 
${\mathcal E}_1\colon{\mathcal P}_{\bf v}\to{\mathcal P}_{\bf u}$ and ${\mathcal E}_2\colon{\mathcal P}_{\bf u}\to{\mathcal P}_{\bf v}$. Hence one also has two second evolute maps
${\mathcal E}_2\circ {\mathcal E}_1: {\mathcal P}_{\bf v}\to{\mathcal P}_{\bf v}$ and ${\mathcal E}_1\circ {\mathcal E}_2: {\mathcal P}_{\bf u}\to{\mathcal P}_{\bf u}$, and these two compositions share the eigenvalues.
\proofend

\subsection{The spectrum of the second evolute transformation}\label{doublespectrum}

Section \ref{iterations} may create an impression that, for a generic spherical polygon $\bf v$, in the case when the maximum module eigenvalue of the transformation ${\mathcal E}^2\colon{\mathcal P}_{\bf v}\to{\mathcal P}_{\bf v}$ is real, the even-numbered evolutes of polygons ${\bf P}\in{\mathcal P}_{\bf v}$ have a prescribed limit shape, that is, there exists a special polygon in ${\mathcal P}_{\bf v}$ such that the even-numbered evolutes of almost any polygon from ${\mathcal P}_{\bf v}$ are asymptotically homothetic to this special polygon. In reality, however, this is not the case: the maximum module eigenvalue of the transformation ${\mathcal E}^2\colon{\mathcal P}_{\bf v}\to{\mathcal P}_{\bf v}$ never has multiplicity one. 

\begin{theorem} \label{spectrumdouble} Let ${\bf v}$ be a generic spherical polygon. Then every non-zero eigenvalue of  the transformation ${\mathcal E}^2\colon{\mathcal P}_{\bf v}\to{\mathcal P}_{\bf v}$ has even multiplicity, generically multiplicity 2. More precisely, the matrix of the restriction of the transformation ${\mathcal E}^2$ to the image of $\mathcal E$ is congugated to the block diagonal matrix of the form $\left[\displaystyle{\begin{array} {cc} A&0\\ 0&A^{\rm T}\end{array}}\right].$\end{theorem}

\proof 
We have two spaces, $V={\mathcal P}_{\bf v}$ and $U={\mathcal P}_{\bf v^\ast}$, and a non-degenerate pairing $\langle\ ,\, \rangle: U \otimes V \to \R$ between them, so we identify $U$ with $V^\ast$. 

 We also have two evolute maps, $A\colon V\to V^\ast$ and $B\colon V^\ast\to V$, which are both anti-self-adjoint: $\langle Ax,x\rangle=0$ and $\langle z,Bz\rangle=0$ for any $x\in V$ and $z\in V^\ast$. Equivalently, the anti-self-adjointness property of $A$ and $B$ may be expressed by the equalities $\langle Ax,y\rangle=-\langle Ay,x\rangle$ and $\langle z,Bw\rangle=-\langle w,Bz\rangle$ for any $x,y\in V$ and $z,w\in V^\ast$.  

Furthermore, we assume that $A$ and $B$ are isomorphisms; it is so if ${\mathcal P}_{\bf v}$ is even-dimensional, and in the odd-dimensional case, we have to replace the spaces ${\mathcal P}_{\bf v}$ and ${\mathcal P}_{\bf v^\ast}$ by the images of the evolute maps.

Let us define a linear symplectic structure on $V$: $\omega(x,y)=\langle Ax,y\rangle$. It is skew-symmetric, because $A$ is anti-self-adjoint, and non-degenerate, because $A$ is an isomorphism. 

A linear transformation $L$ of a symplectic space is called skew-Hamiltonian, if $\omega(Lx,y)=\omega(x,Ly)$  for all $x,y$ (compare with a more familiar definition of a Hamiltonian map $H$, an infinitesimal symplectomorphism, satisfying the condition $\omega(Hx,y)=-\omega(x,Hy)$).

We claim that the map $BA\colon V\to V$ is skew-Hamiltonian. Indeed, 
$$
\omega(x,BAy)=\langle Ax,BAy\rangle=-\langle ABAy,x\rangle=\langle ABAx,y\rangle=\omega(BAx,y)
$$
 (the first and last equalities hold by the defiinition of $\omega$, and the second and third equalities are the anti-self-adjointness properties of $A$ and $ABA$). 

Now the last and the most important step. Theorem 1 of  \cite{FMMX} asserts that the matrix of a skew-Hamiltonian transformation is conjugated by a symplectic matrix to a matrix of the form stated in Theorem \ref{spectrumdouble}. \proofend

\begin{remark}
{\rm In dimension two, the inscribed polygons ${\bf v}$ and its dual ${\bf v^\ast}$ differ by rotation through $90^{\circ}$, and the evolute map depends only on the exterior angles of a polygon. This makes it possible to identify the space ${\mathcal P}_{\bf v}$ with ${\mathcal P}_{\bf v^\ast}$, and the evolute map $A$ with $B$. Then the matrix of the second evolute map becomes the square of the matrix of the evolute map. Theorem 7 of \cite{ArFuIzTaTs} states that the spectrum of this evolute map is symmetric with respect to the origin and, furthermore, the opposite eigenvalues have the same geometric multiplicity and the same sizes of the Jordan blocks. This implies the 2-dimensional case of the above Theorem \ref{spectrumdouble}.
}
\end{remark}

\begin{remark}
{\rm Theorem \ref{spectrumdouble} implies Theorems \ref{pentatheorem} and \ref{hexatheorem}. Indeed, in these two theorems one deals with maps of two-dimensional spaces, and Theorem \ref{spectrumdouble} implies that these maps are diagonal.
}
\end{remark}

\section{Curves}\label{smooth}

In this section, we restrict ourselves to the three-dimensional case. We hope to consider the situation in higher dimensions, as well as in the elliptic and hyperbolic geometries, in the near future.

Recall that the evolute of a spacial curve is the locus of the centers of its osculating spheres or, equivalently, the enveloping curve of the 1-parameter family of its normal planes.

\subsection{Spherical curves}\label{spherical} 
Let us recall some facts about spherical curves; see, e.g., \cite{Arn}. 

By a {\it spherical curve} we always mean a closed immersed locally convex smooth curve in the sphere $S^2$; ``smooth" here and below means ${\mathcal C}^\infty$; ``locally convex" or ``inflection free" means that the geodesic curvature of the curve does not vanish, that is, the curve is nowhere abnormally well approximated by a great circle. 

A coorientation of a spherical curve is a choice of a unit normal vector field. An oriented spherical curve has a canonical coorientation obtained by rotating the orienting vector $90^{\circ}$ in the positive direction.

Let $\gamma$ be a cooriented spherical curve. The dual curve $\overline\gamma$ is obtained by moving each point of $\gamma$ distance $\pi/2$ along the great circle in the direction of the coorientation. For example, the Southward cooriented Arctic Circle is dual to the Tropic of Capricorn, and the Northward cooriented Arctic Circle is dual to the Tropic of Cancer. 

The class of immersed locally convex curves is invariant with respect to duality. Moreover, for a cooriented spherical curve $\gamma$, the curve $\overline\gamma$ possesses a canonical coorientation, and the second dual curve is antipodal to the original one. 

Let $\gamma(t)$ be an arc length parameterized spherical curve with the geodesic curvature $\kappa$, let $\overline\gamma$ be its dual curve, and $s$ the arc length parameter on $\overline\gamma$.

\begin{lemma} \label{stot}
One has:
$$
\frac{ds}{dt} = |\kappa|.
$$
\end{lemma} 

\proof
The coorientation of the curve $\gamma$ is given by the unit normal vector $\gamma \times \gamma_t$, hence $\overline\gamma=\gamma \times \gamma_t$. Denote the $90^{\circ}$ rotation of tangent vectors in the positive direction by $J$. Then one has:
$$
\overline\gamma_t = \gamma \times \gamma_{tt} = \kappa\ \gamma \times J(\gamma_{t}) = -\kappa \gamma_t,
$$
and hence $ds/dt = |\overline\gamma_t|=|\kappa|$.
\proofend

We will need the following property of the tangent indicatrix of a spherical curve: its total geodesic curvature is zero (equivalently, the appropriately defined, area bounded by the tangent indicatrix of a spherical curve is a multiple of $2\pi$), see \cite{Arn,So}. 

\subsection{Spacial curves with a fixed tangent indicatrix} \label{fixdtan}

Let $\gamma(t)$ be an arc length parameterized spherical curve.  Denote the unit position and velocity vectors $\gamma$ and $\gamma_t$ by $T$ and $N$, respectively, and let $B = T \times N$. We wish to describe the spacial curves for which $\gamma$ serves as the tangent indicatrix.

Let $\rho(t)$ be a smooth periodic function with the period equal to the length of $\gamma$, satisfying the relation 
\begin{equation} \label{closed}
\int\rho(t)\gamma(t)\ dt = 0.
\end{equation}
Define a spacial curve by the formula
$$\Gamma(t)=\Gamma(0)+\displaystyle\int_0^t\rho(s)\gamma(s)\, ds.$$
Due to (\ref{closed}), $\Gamma$ is closed, and it is well defined up to a parallel translation. The curve $\Gamma$ receives an orientation from the curve $\gamma$.

The curve $\Gamma$ consists of alternating arcs, that we call {\it positive} and {\it negative}, according to the sign of the function $\rho$.  The vector 
$\rho(t)\gamma(t)$ is the  velocity vector of the curve $\Gamma(t)$ that ``flips" at zeros of $\rho$. 

Assume that the function $\rho$ has only simple zeroes. The corresponding points of the curve $\Gamma$ are simple cusps,  
diffeomorphic to the curve $(u^2,u^3,u^4)$ that has a cusp at the origin. One can choose  local coordinates $X,Y,Z$ in space so that the curve is given by the parametric equations 
\begin{equation} \label{cusps}
X(u)=au^2 +O(u^3),\, Y(u)=bu^3 + O(u^4),\, Z(u)=cu^4 + O(u^5),
\end{equation}
where $a,b,c$ are non-zero real numbers. 

The curve $\Gamma$ is equipped with an orthonormal frame $(T,N,B)$; the vectors of this frame agree at the cusps. The vector $T$ defines a map from $\Gamma$ to the unit sphere and, by construction, the tangent indicatrix of $\Gamma$ is the spherical curve $\gamma$. Since the spherical curve $\gamma$ is smooth and locally convex, the spacial curve $\Gamma$ has non-vanishing curvature and torsion.

We call the spacial curves, resulting from the above described construction with $\rho$ having only simple zeros, {\it admissible curves}.  
To recap, an admissible curve is oriented, has a finite number of simple cusps, non-vanishing curvature and torsion, and its pieces between the cusps are marked as positive and negative. An admissible curve is equipped with an orthonormal frame $(T,N,B)$, and the unit vector $T$ agrees with the orientation on the positive pieces, and is opposite to it on the negative ones. 

Denote by ${\mathcal C}_\gamma$  the space of  smooth periodic function $\rho(t)$ satisfying relation (\ref{closed}).  This space is a continuous analog of the space ${\cal P}_{\bf v}$. The admissible curves having  $\gamma$ as the tangent indicatrix play the role of non-degenerate polygons in ${\mathcal P}_{\bf v}$; they form an open dense subset of the space ${\mathcal C}_\gamma$. Of course,  spacial curves  corresponding to some $\rho\in {\mathcal C}_\gamma$ may be severely degenerate, even consisting of one point (for $\rho=0$). 

Sometimes we shall abuse notation and write $\Gamma \in {\mathcal C}_\gamma$ to indicate that the spacial curve $\Gamma$ is the result of the above described construction.

\subsection{Differential geometry of admissible curves} \label{diffgeo}

Let $\gamma(t)$ be an arc length parameterized  spherical curve, $\rho(t)$ a smooth function on $\gamma$ with only simple zeroes, and let $\Gamma$ be the respective admissible curve. Let $x$ be the arc length parameter on $\Gamma$, and $k>0$ and $\tau$ be its curvature and torsion. Let $\sigma(x)=\pm 1$ be a locally constant function, equal to 1 on the positive, and to $-1$ on the negative pieces of  $\Gamma$.

The next lemma is proved by a direct calculation which we omit.

\begin{lemma} \label{atcusp}
The curvature $k(u)$ and the torsion $\tau(u)$ at the cusp (\ref{cusps}) are both infinite: 
$$
k(u)\sim\displaystyle\frac{3b}{4a^2|u|},\, \tau(u)\sim\displaystyle\frac{4c}{3abu}.
$$
\end{lemma}

In particular, the torsion changes sign at cusps.

\begin{lemma} \label{frame}
The orthonormal frame $(\sigma T, \sigma N, B)$ is the Frenet frame along $\Gamma$, with the
Frenet equations taking the form
\begin{equation} \label{newFr}
T_x = \frac{\sigma}{\rho} N,\ N_x = -\frac{\sigma}{\rho} T + \sigma \tau B,\ B_x = - \sigma \tau N.
\end{equation}
One has:
$$
\frac{dx}{dt}=\frac{1}{k} = \sigma \rho,\, \frac{dt}{dx}=k = \frac{\sigma}{\rho}.
$$
\end{lemma}

\proof On a positive piece of the curve, $T = \Gamma_x$, and hence $(T, N, B)$ is the Frenet frame. At a cusp, the vectors $T$ and $N$ remain the same, but the velocity and acceleration vectors of $\Gamma$ change signs; the binormal vector is given by the cross-product, and it remains the same. Rewriting the usual Frenet equations for the Frenet frame $(\sigma T, \sigma N, B)$ yields (\ref{newFr}). 

One has $\Gamma_t = \rho \gamma$ by construction, hence $x_t = |\rho|$. On the other hand, $T_x =kN$, and 
$1 = |T_t| = |T_x| x_t = k x_t,$ hence $1/k=|\rho|$. One takes care of the sign by multiplying $\rho$ by $\sigma$.
\proofend

 It follows from the previous two lemmas that, at a cusp, the product $\rho(u)\tau(u)$ has a finite non-zero limit. We can say more.

\begin{lemma} \label{geodcurv}
The geodesic curvature $\kappa$ of the spherical curve $\gamma$ equals $\rho \tau$.
\end{lemma}

\proof
For an arc length parameterized spherical curve $\gamma(t)=T$, the vector of the geodesic curvature $\kappa B$ equals $\gamma+\gamma_{tt}=T+N_t$. Using Lemma \ref{frame}, one finds
$$
N_t =   N_x x_t =  \left( -\frac{\sigma}{\rho} T + \sigma \tau B\right) \sigma \rho = -T + \rho \tau B.
$$
Hence $\kappa = \rho \tau$.
\proofend

\subsection{Evolutes}\label{evolutesmooth}

Let $\Gamma \in {\mathcal C}_{\gamma}$ be an admissible curve and let $\overline\Gamma$ be the evolute of $\Gamma$.

\begin{lemma} \label{evolform}
One has
$$
\overline\Gamma = \Gamma + \rho N + \sigma \frac{\rho_x}\tau B,
$$
and
$$
\overline\Gamma_x = \sigma \left[ \rho\tau+\left(\frac{\rho_x}\tau\right)_x \right] B.
$$
\end{lemma} 

\proof
The formula for the evolute of a curve in terms of its Frenet frame can be found in, e.g., \cite{St,UrVa} and \cite{Fu} (this formula was originally due to Monge):
$$
\overline\Gamma = \Gamma + \frac{1}{k} N + \frac{({1}/{k})_x}{\tau} B.
$$
We need to adjust this formula taking Lemma \ref{frame} into account, that is, by replacing $N$ with $\sigma N$ and $1/k$ with $\sigma \rho$.
This yields the first formula. The second one is obtained by differentiation, using formulas (\ref{newFr}).
\proofend

\begin{remark} \label{atcusprmk}
{\rm At a cusp $(at^2,bt^3,ct^4)$,  the quantity $\rho_x/\tau$ has a finite limit, and the center of the osculating sphere is located on the binormal, the vertical axis, at distance $a^2/(2c)$ from the origin. The cusps of $\Gamma$ correspond to regular points of the evolute $\overline\Gamma$.
}
\end{remark}

Let $\overline\gamma$ be the spherical curve dual to $\gamma$ and $s$ its arc length parameter. 
Since the binormal vector $B$ is the position vector of the curve  $\overline\gamma$, the second formula of Lemma \ref{evolform} implies 

\begin{corollary} \label{dualind}
The tangent indicatrix of the evolute is spherically dual to the tangent indicatrix of the original curve.
\end{corollary}

\begin{remark}
{\rm
The second formula of Lemma \ref{evolform} also implies that 
$$
\Phi:= \rho\tau+\left(\frac{\rho_x}\tau\right)_x =0
$$
is a criterion  for a spacial curve to lie on a sphere. This is a classical result, see, e.g., \cite{St}.
}
\end{remark}

We want to find  a function $\overline\rho(s) \in {\mathcal C}_{\overline\gamma}$ that yields the evolute $\overline\Gamma$, as described in Section \ref{fixdtan}. Let us assume that the geodesic curvature $\kappa$ of the curve $\gamma$ is positive (if not,   reverse the orientation of the curve). 

\begin{theorem}\label{linmap} 
One has:
$$
\overline\rho = \rho + \frac1\tau\left(\frac{\rho_x}\tau\right)_x = \rho + \rho_{ss}.
$$
\end{theorem}

\proof
The function $\overline\rho$ is characterized by the  equality $\overline\Gamma_s = \overline\rho(s) {\overline\gamma} (s)$. 

First, we compute
\begin{equation} \label{xtos}
x_s = x_t\ t_s = (\sigma \rho)\ \frac{1}{\kappa} = \sigma \rho\ \frac{1}{\rho \tau} = \frac{\sigma}{\tau},
\end{equation} 
where the second equality makes use of Lemmas  \ref{frame} and \ref{stot}, and the third of Lemma \ref{geodcurv}. 
Next,
$$
\overline\Gamma_s = \overline\Gamma_x\ x_s = (\sigma \Phi B)\ \frac{\sigma}{\tau} =  \frac\Phi\tau\ B,
$$
implying the first equality of the theorem. Equation (\ref{xtos}) implies that
$$
\frac{d}{ds} = \frac{\sigma}{\tau}\ \frac{d}{dx}.
$$
Applying this differential operator twice to $\rho$ yields the second equality of the theorem.
\proofend

Theorem \ref{linmap} describes the evolute map on admissible curves as a linear transformation $\rho \mapsto \rho + \rho_{ss}$. This formula extends to the whole space ${\mathcal C}_{\gamma}$ and yields the linear evolute map ${\mathcal E}\colon C_\gamma\to C_{\overline\gamma}$. This result may be regarded as a continuous analog of Theorem \ref{linear}.

\subsection{The kernel and the image of the evolute map} \label{kerim}

This section is a continuous counterpart to Section \ref{rank}. 

Let $\gamma$ be a generic  spherical curve. Consider the evolute map ${\cal E}: {\cal C}_{\gamma} \to {\cal C}_{\overline\gamma}$. The next proposition is an analog of Proposition \ref{kernel}: the infinite-dimensional space ${\cal C}_{\gamma}$, in a sense, is odd-dimensional.  

\begin{proposition} \label{infodd}
The map ${\cal E}$ is a linear bijection.
\end{proposition}

\proof
Similarly to the polygonal case, $\Ker {\cal E}$ consists of the  curves in ${\cal C}_{\gamma}$ that lie on a sphere. As we mentioned earlier, the tangent indicatrix of a spherical curve has zero total geodesic curvature. This condition fails for a generic $\gamma$, hence  $\Ker {\cal E} = 0$.

Given $\overline\Gamma \in {\cal C}_{\overline\gamma}$, we want to construct its involute, a curve $\Gamma$ whose evolute is $\overline\Gamma$. As we mentioned in Introduction, the evolute of a curve is the envelope of its normal planes. That is, $\Gamma$ is normal to the family of the osculating planes of $\overline\Gamma$.

Let $y$ be a parameter on $\overline\Gamma$, and consider the family $\xi(y)$ of the osculating planes of $\overline\Gamma$. Since ${\overline\gamma}$ is a closed curve, the curve $\overline\Gamma$ has an even number of cusps. The  frame $(T,N)$ gives the osculating planes of  $\overline\Gamma$ consistent orientations. Thus $\xi(y)$ is a loop in the space of oriented planes in $\R^3$, and we need to find a closed curve, normal to this 1-parameter family of planes.

Start with some plane $\xi(y_0)$, and pick a point $A$ in this plane. There is a unique curve through $A$, orthogonal to our family of planes. After one traverses the curve $\overline\Gamma$, the orthogonal curve returns to the same plane $\xi(y_0)$, say, at point $B$. This defines a map $A \mapsto B$ of this plane, and we want to find its fixed point.

We claim that this map of the  plane $\xi(y_0)$ is an orientation-preserving isometry. It suffices to establish an infinitesimal version of this claim. Indeed, the map $\xi(y) \to \xi(y+dy)$, given by the orthogonal curves to the family of planes, is a rotation about the intersection line of these two infinitesimally close planes (this line is tangent to the curve $\overline\Gamma$).

Thus we have an orientation preserving isometry of the plane $\xi(y_0)$, generically, a rotation. A rotation has a unique fixed point, as needed.
\proofend

\begin{remark} 
{\rm One can show that the angle of the rotation $\xi(y_0) \to \xi(y_0)$ equals the total curvature of the curve $\overline\Gamma$.
}
\end{remark}

\subsection{${\mathcal C}_\gamma$ as dual space to ${\mathcal C}_{\overline\gamma}$}\label{dualsectcurves}

This section is a continuous counterpart to Section \ref{dualsect}: for a generic spherical curve $\gamma$, we will define a non-degenerate pairing between ${\mathcal C}_\gamma$ and ${\mathcal C}_{\overline\gamma}$ and will prove that the evolute transformation in anti-self-adjoint. 

We continue using the same notations as in the preceding subsections.

The {\it support number} $\lambda(x)$ is defined as the signed distance from the origin to the osculating plane of $\Gamma$ at the point $\Gamma(x)$, specifically, $\lambda(x)=\Gamma(x)\cdot B(x)$. For $\Gamma\in{\mathcal C}_\gamma, \overline\Gamma\in{\mathcal C}_{\overline\gamma}$, we put
\begin{equation} \label{scal}
\langle\Gamma,\overline\Gamma\rangle=\int_{\overline\gamma}\lambda(x(s))\overline\rho(s)\, ds.
\end{equation}
(Notice that $\rho=dx/dt$ is an analog of the side lengths for polygons.)

The next proposition is a continuous analog of Lemma \ref{welldef}, and Propositions \ref{antisym} and \ref{nondeg}.

\begin{proposition}\label{nondegcurves} For a generic $\gamma$, the pairing $\langle\ \rangle$ is well-defined, anti-symm\-etric, and non-degenerate.\end{proposition}

\proof 
We claim that 
$$
\int_{\overline\gamma}\lambda(x(s))\overline\rho(s)\, ds =
\int_{\overline\gamma} \Gamma \cdot \overline\Gamma_s\, ds.
$$
Indeed, 
$\overline\Gamma_s =  \overline\rho B$, 
hence $\Gamma \cdot \overline\Gamma_s = (\Gamma \cdot B) \overline\rho = \lambda \overline\rho.$

It follows that a parallel translation through vector $R$ results in the following change of (\ref{scal}):
$$
\int_{\overline\gamma} R \cdot \overline\Gamma_s\, ds = R \cdot \int_{\overline\gamma} \overline\Gamma_s\, ds = R \cdot \int \overline\rho\ \overline\gamma\, ds = 0.
$$
Therefore the pairing does not depend on the choice of the origin.

Concerning anti-symmetry,
$$
\langle\Gamma,\overline\Gamma\rangle= \int \Gamma \cdot \overline\Gamma_s\, ds = \int  \Gamma \cdot d\overline\Gamma = - \int \overline\Gamma \cdot d\Gamma = - \langle\overline\Gamma,\Gamma\rangle.
$$

Next, we prove non-degeneracy. 
Given $\Gamma\in {\mathcal C}_{\gamma}$, we want to find  $\overline\Gamma \in {\mathcal C}_{\overline\gamma}$ such that $\langle\Gamma,\overline\Gamma\rangle\neq0$; in other words, we want to find a function $\overline\rho(s)$  such that 
$$
\int_{\overline\gamma}\overline\rho\, \overline\gamma\, ds=0,\ {\rm but}\ \int_{\overline\gamma}\lambda\, \overline\rho\, ds\neq0.
$$

It is clear that we can find such  $\overline\rho$ if and only is the function $\lambda(s)$ is not a linear combination with constant coefficients of the three components of the vector function $\overline\gamma$, that is, the binormal vector $B$.

The function $\Gamma \cdot B$ is a linear combination of the components of the vector $B$ if and only if $\Gamma \cdot B = C \cdot B$ for some constant vector $C$, that is, when $(\Gamma - C) \cdot B=0$, or $(\Gamma - C) \cdot (T\times N)$, that is,  $\det (\Gamma - C, \Gamma_x, \Gamma_{xx}) = 0$ identically.

The last identity implies that the curve $\Gamma$ is planar, which is not the case if $\gamma$ is generic.
 \proofend

Assume now that $\overline\Gamma$ is the evolute of $\Gamma$. 
The next result is a continuous analog of Proposition \ref{orthogonal}. 

\begin{proposition}\label{anrtiselfadjcurves} 
One has $\langle\Gamma,\overline\Gamma\rangle=0$.
\end{proposition}

\proof One has
\begin{equation*} 
\langle\Gamma,\overline\Gamma\rangle =  \int (\Gamma\cdot B)\, (\rho + \rho_{ss})\, ds =
\int [(\Gamma\cdot B) + (\Gamma_{ss}\cdot B) + 2(\Gamma_s\cdot B_s) + (\Gamma\cdot B_{ss})]\, \rho\, ds,
\end{equation*}
where the last equality is the result of integration by parts twice.

Let us examine the integrands. One has
$$
\Gamma_s = \Gamma_t\, t_s = \frac{1}{\rho \tau}\, \rho T = \frac{1}{\tau}\, T
$$
and, since $N_t = -T + \kappa B$,
$$
B_s = (T \times N)_s = T \times N_t\ t_s =  \frac{1}{\kappa}\ \kappa T \times B = -N.
$$
It follows that $\Gamma_s \perp B_s$, and the third integrand vanishes. We continue:
$$
\Gamma_{ss} = \left(\frac{1}{\tau}\, T\right)_t\, t_s,
$$
which is a linear combination of $T$ and $N$, hence orthogonal to $B$. Thus the second integrand vanishes as well.
Next,
$$
B_{ss} = - N_t\, t_s = \frac{1}{\kappa}\, (T - \kappa B) = \frac{1}{\kappa}\, T - B,\ \ {\rm hence}\ \ 
(\Gamma\cdot B) + (\Gamma\cdot B_{ss}) = \frac{1}{\kappa}\, T.
$$
Thus what remains of the integral is
$$
\int \frac{1}{\kappa}\, (\Gamma\cdot T)\, \rho\, ds = \int (\Gamma \cdot \Gamma_t)\, dt = \frac{1}{2} \int (\Gamma \cdot \Gamma)_t\, dt = 0,
$$
where the first equality is due to the facts that $\Gamma_t = \rho T$ and $ds/\kappa = dt$.
\proofend

\subsection{Spacial hypocycloids} \label{spathypo}

Consider an arc length parameterized circle of latitude $\gamma$ at height $\sqrt{1-r^2}$:
\begin{equation} \label{circle}
x=r\cos\left(\frac tr\right),\, y=r\sin\left(\frac tr\right),\ z=\sqrt{1-r^2},
\end{equation}
with $0\leq t\leq2\pi r.$ Let $\rho(t)\in {\mathcal C}_{\gamma}$ be a $2\pi r$-periodic function. Then
$$
\int_0^{2\pi r}\cos\left(\frac tr\right)\, \rho(t)\, dt=\int_0^{2\pi r}\sin\left(\frac tr\right)\, \rho(t)\, dt=\int_0^{2\pi r}\rho(t)dt=0.
$$
Set $\alpha=\displaystyle\frac tr;$ then
$$
\int_0^{2\pi}\cos\alpha\ \rho(\alpha)\ d\alpha =\int_0^{2\pi}\sin\alpha\ \rho(\alpha)\ d\alpha=\int_0^{2\pi}\rho(\alpha)\ d\alpha=0,
$$
that is, the Fourier expansion of $\rho(\alpha)$ is free from the constant term and from the first harmonics.

We mention, in passing, that this implies that the function $\rho(\alpha)$ has at least four zeros on the interval $[0, 2\pi)$; see, e.g., \cite{FT}. 

Next, consider the arc length parameter $s$ on the dual curve $\overline\gamma$, a
circle of latitude at height $r$. One has
$$
s=\frac{\sqrt{1-r^2}}r=\sqrt{1-r^2}\ \alpha,
$$
hence
\begin{equation} \label{Lap}
\overline\rho=\rho+\rho_{ss}=\rho+\frac1{1-r^2}\rho_{\alpha\alpha}.
\end{equation}

Consider now the case when the function $\rho$ is a pure harmonic: $\rho(\alpha)=\cos k\alpha$ with integral $k\geq2$. Then, by formula (\ref{Lap}), $\overline\rho$ is proportional to $\rho$. In this case, the parametric equations of the curve $\Gamma(t)$ (with the indicatrix being the circle of latitude (\ref{circle})) are
\begin{equation} \label{hypo}
\begin{split}
& x=r\left(\frac{\sin(k-1)t}{k-1}+\frac{\sin(k+1)t}{k+1}\right),\\
& y=r\left(\frac{\cos(k-1)t}{k-1}-\frac{\cos(k+1)t}{k+1},\right),\ z=2\sqrt{1-r^2}\ \frac{\sin kt}{k}.
\end{split}
\end{equation}
It follows from formula (\ref{Lap}) (and can be confirmed by a direct computation) that the evolute $\overline\Gamma$ is obtained from $\Gamma$ by switching $r\longleftrightarrow\sqrt{1-r^2}$ and a homothety with the coefficient $\displaystyle\frac{r^2+k^2-1}{r\sqrt{1-r^2}}.$ 

\begin{corollary} \label{selfdual}
The second evolute $\overline{\overline\Gamma}(t)$ is obtained from $\Gamma(t)$ by a  homothety with the coefficient $\displaystyle\frac{r^2(1-r^2)+k^2(k^2-1)}{r^2(1-r^2)}$. 
\end{corollary}

This makes the curve $\Gamma(t)$ similar to the classical hypocycloid.

\begin{figure}[hbtp]
\centering
\includegraphics[height=3in]{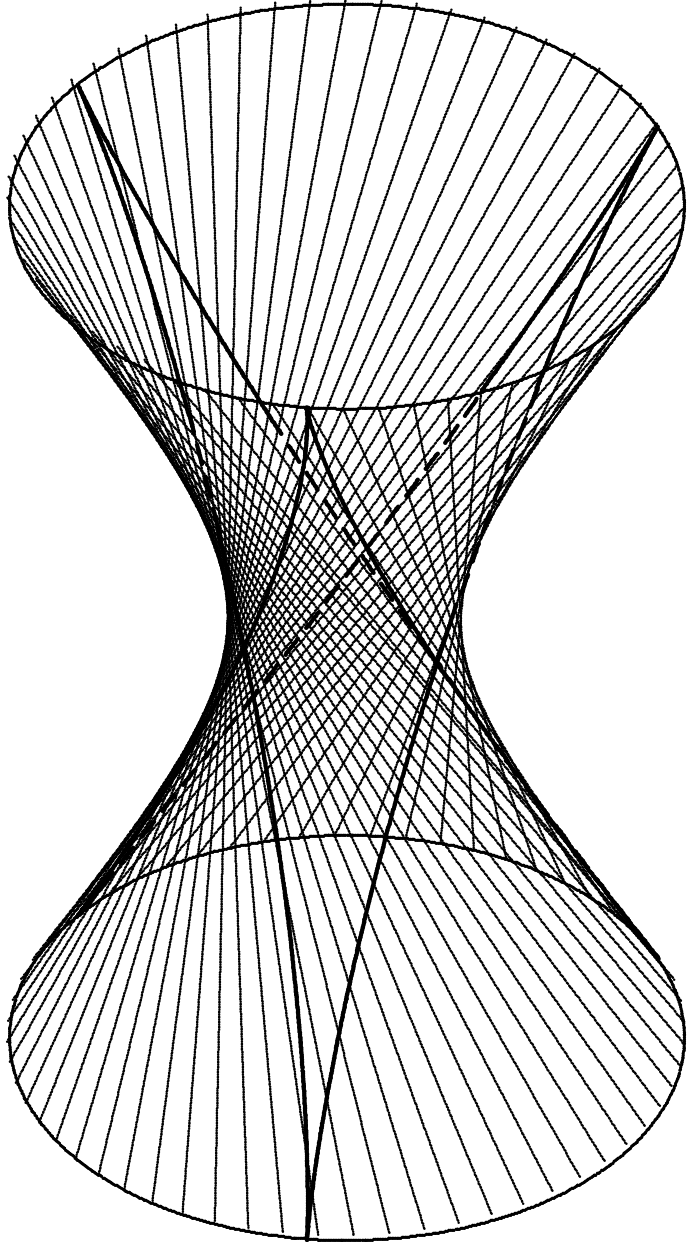}\qquad\qquad\qquad
\includegraphics[height=3in]{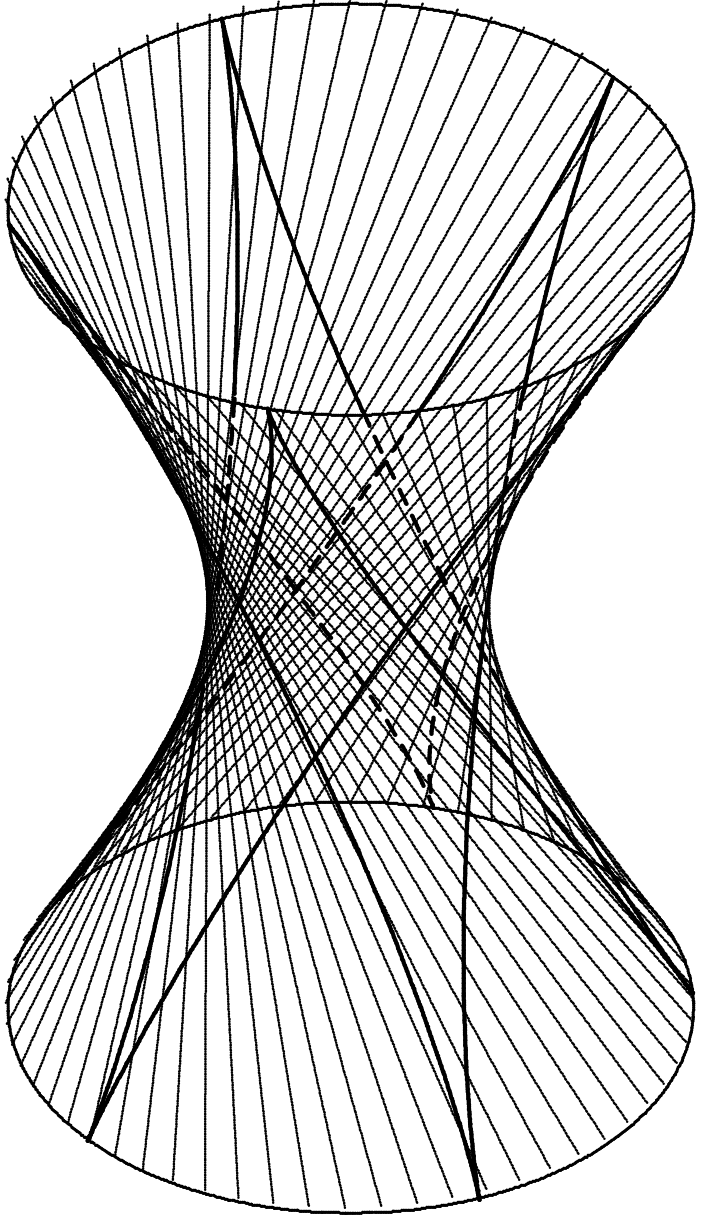}
\caption{spacial hypocycloids: $k=3$ and $k=\displaystyle\frac52$}
\label{hypocycloids}
\end{figure}

Let us provide some geometric information regarding these ``hypocycloids." The next proposition is proved by a straightforward calculation.

\begin{proposition} \label{hyperb}
The curve $\Gamma(t)$ with the parametric equations (\ref{hypo}) is contained in the hyperboloid $$\frac{k^2-1}{4r^2}(x^2+y^2)-\frac{k^2}{4(1-r^2)}\,z^2=\frac1{k^2-1}$$between the planes $z=\pm\displaystyle\frac{2\sqrt{1-r^2}}k$. It has $2k$ cusps at the points$$\left(2kr\cos\frac{2i-1}{2k},\, 2kr\sin\frac{2i-1}{2k},\, 2(-1)^i\frac{\sqrt{1-r^2}}k\right),\ i=1,2,\dots,2k.$$
\end{proposition}

The number $k$ may be rational, $k=p/q>1,\, (p,q)=1.$ In this case, the indicatrix will be the circle (\ref{circle}) traversed $q$ times, the curve $\Gamma$ will be contained in the same hyperboloid and will have $2p$ cusps. See Figure \ref{hypocycloids}.

\end{document}